\documentclass{amsart}
\usepackage{amssymb}
\usepackage{amsthm}
\usepackage[dvips]{graphicx}

\theoremstyle{definition} \newtheorem{defin}{Definition}[section]
\newtheorem{ex}{Example}[section]
\theoremstyle{remark} \newtheorem{rem}{Remark}[section]
\newcommand{\pn}{\par\noindent} \newcommand{\pmn}{\par\medskip\noindent}

\begin{document}
\title{A two dimensional analog of the three gap theorem}
\author{Yury Kochetkov \and Alexandr Osipov}
\date{}
\begin{abstract} For a two dimensional vector $\bar v=(\alpha,\beta)$, where
$\alpha>0$, $\beta>0$ are irrational numbers, independent over
$\mathbb{Q}$, we consider the set $D_n=\{(i\alpha\, {\rm mod }\,
1,i\beta\, {\rm mod}\, 1),\,\,i=1,\ldots,n\}$ in a two dimensional torus
and the partition of this torus into Voronoi cells. Areas and forms of
these cells are the subject of this experimental work.
\end{abstract}
\email{yukochetkov@hse.ru, afosipov@edu.hse.ru} \maketitle

\section{Introduction}
\pn Three gap theorem (Steinhaus conjecture) was proposed by
Polish mathematician Hugo Steinhaus in the 30th years of the 20th
century. Let $\alpha>0$ be an irrational number. Let us consider
the sequence $A_n=\{\alpha,2\alpha,\ldots,n\alpha\}$ and the
sequence of fractional parts
$B_n=\{\{\alpha\},\{2\alpha\},\ldots,\{n\alpha\}\}\subset [0,1]$.
We identify the end points of the segment $[0,1]$ and get the
circle $S^1$ of the unit length. Points of the set $B_n$ divide
$S^1$ into $n$ segments. The Steinhaus conjecture states that
these segments have not more than 3 different lengths. The
Steinhaus conjecture was proved in 1958 by Hungarian mathematician
V. S\'{o}s \cite{Sos} and was reproved many times (see \cite{Rav},
\cite{MS}, for example). \pmn In this purely experimental work we
consider a two-dimensional analog of the Steinhaus conjecture,
where instead of lengths of intervals we consider areas of Voronoi
cells. \pn
\begin{defin} Let $A=\{a_1,\ldots,a_n\}$ be a finite set of points in the
plane. The Voronoi cell $V_i$ of an element $a_i$ is defined in
the following way:
$$V_i=\{x\in\mathbb{R}^2\,|\,|x-a_i|\leqslant |x-a_j|\,\,\forall j\neq i\}.$$
Voronoi cell is a convex polygon, maybe unbounded. \end{defin} \pn
\begin{rem} As our ground set will be a two dimensional torus, then
our Voronoi cells will be bounded convex polygons. \end{rem} \pn
Let $\alpha$ and $\beta$ be positive, irrational numbers,
independent over $\mathbb{Q}$. Let us consider the sequence
$C_n=\{\bar v,2\bar v,\ldots,n\bar v\}\subset \mathbb{R}^2$, where
$\bar v=(\alpha,\beta)$ --- a two dimensional vector, and the set
$D_n=\{({\alpha},{\beta}),(\{2\alpha\},\{2\beta\}),\ldots,(\{n\alpha\},
\{n\beta\})\}\subset [0,1]\times [0,1]$. Let us identify left and
right sides of square and upper and lower sides. We obtain a torus
and set of points in it. The Voronoi cells give a partition of the
torus into disjoint convex polygons. \pn
\begin{ex} Here is a partition of torus into Voronoi cells for
the vector $v=(\sqrt 2,\sqrt 3)$ and $n=20$:
\begin{center}
\includegraphics[scale=0.5]{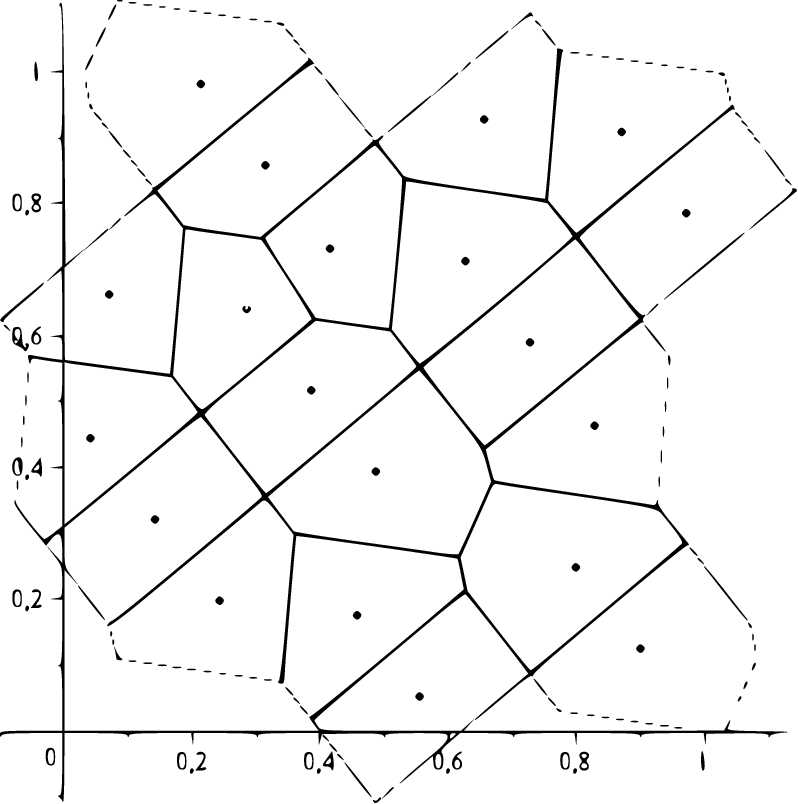}
\end{center}
\begin{center}{Figure 1}\end{center}
\end{ex}
\pn Let $S(n)$ be the number of \emph{different} areas of Voronoi
cells of the set $D(n)$. The main question: is the sequence $S(n)$
bounded or not? We computed areas of Voronoi cells with $10^{-80}$
precision for $n\leqslant 1000$. Results are presented in tables
below. By $M_k(n)$ we denote the number of $k$-polygons.

\section{Results}
\subsection{$v=(\sqrt 2,\sqrt 3)$}
\[\begin{tabular}{|c|c|c|c|c|c|c|c|c|c|c|c|c|c|c|c|c|c|}\hline
n&20&30&40&50&60&70&80&90&100&110&120&130&140&150&200&500&1000\\
\hline $S(n)$&6&6&6&6& 6&8&8&7&8& 7&7&7&7&7& 7&7&10\\ \hline
$M_3(n)$&0&0&0&0& 0&0&0&0&0& 0&0&0&0&0& 0&104&0\\ \hline
$M_4(n)$&0&0&0&0& 0&0&0&8&0& 16&6&0&0&0& 0&38&0\\ \hline
$M_5(n)$&6&8&12&30& 24&24&30&14&18& 12&32&44&44&44& 40&276&328\\
\hline
$M_6(n)$&8&14&18&2& 20&22&20&44&0& 38&38&42&52&62& 120&0&508\\
\hline $M_7(n)$&6&8&8&6& 8&24&30&18&46& 44&44&44&44&44& 40&0&82\\
\hline $M_8(n)$&0&0&2&12& 8&0&0&6&36& 0&0&0&0&0& 0&82&0\\ \hline
$M_9(n)$&0&0&0&0& 0&0&0&0&0& 0&0&0&0&0& 0&0&82\\ \hline
\end{tabular}\] \pmn \begin{center}{Table 1}\end{center}
\subsection{$v=(\sqrt 2,\sqrt 5)$}
\[\begin{tabular}{|c|c|c|c|c| c|c|c|c|c| c|c|c|c|c| c|c|c|}\hline
n&20&30&40&50& 60&70&80&90&100& 110&120&130&140&150& 200&500&1000\\
\hline $S(n)$&6&7&7&7& 6&8&8&8&8& 7&9&7&6&8& 8&6&6\\ \hline
$M_4(n)$&0&0&0&16& 26&22&12&2&0& 0&0&0&0&0& 0&0&0\\ \hline
$M_5(n)$&8&8&6&0& 0&14&34&54&46& 34&34&34&24&24& 62&102&170\\
\hline $M_6(n)$&4&18&0&2& 0&0&0&0&20& 50&62&72&92&102&
76&296&660\\ \hline $M_7(n)$&8&0&22&32& 16&10&10&10&22&
18&14&14&24&24& 62&102&170\\ \hline $M_8(n)$&0&4&12&0&
18&24&24&24&12& 8&10&10&0&0& 0&0&0\\ \hline
\end{tabular}\] \pmn \begin{center}{Table 2}\end{center}
\subsection{$v=(\sqrt 3,\sqrt 5)$}
\[\begin{tabular}{|c|c|c|c|c|c|c|c|c|c|c|c|c|c|c|c|c|c|}\hline
n&20&30&40&50&60&70&80&90&100&110&120&130&140&150&200&500&1000\\
\hline $S(n)$&5&5&9&8& 5&8&8&8&8& 8&8&6&6&6& 7&6&7\\ \hline
$M_4(n)$&0&0&8&8& 0&0&0&0&0& 0&0&0&0&0& 0&0&246\\ \hline
$M_5(n)$&2&8&8&4& 6&14&16&14&14& 16&14&6&26&46& 60&60&0\\ \hline
$M_6(n)$&16&14&8&28& 48&42&48&68&72& 78&92&118&88&58& 86&380&262\\
\hline $M_7(n)$&2&8&8&0& 6&14&16&2&14& 16&14&6&26&46& 48&60&492\\
\hline $M_8(n)$&0&0&8&10& 0&0&0&6&0& 0&0&0&0&0& 6&0&0\\ \hline
\end{tabular}\] \pmn \begin{center}{Table 3}\end{center}
\subsection{$v=(\sqrt 5,\sqrt 6)$}
\[\begin{tabular}{|c|c|c|c|c|c|c|c|c|c|c|c|c|c|c|c|c|c|}\hline
n&20&30&40&50&60&70&80&90&100&110&120&130&140&150&200&500&1000\\
\hline $S(n)$&5&6&6&7& 6&7&5&6&6& 6&6&6&5&6& 6&6&7\\ \hline
$M_4(n)$&2&0&0&0& 0&0&4&0&0& 0&0&0&0&0& 0&0&0\\ \hline
$M_5(n)$&0&2&4&18& 18&12&0&2&22& 42&62&76&76&76& 76&76&178\\
\hline $M_6(n)$&14&26&32&14& 26&52&68&86&56& 34&14&6&26&26&
48&348&644\\ \hline $M_7(n)$&4&2&4&18& 14&0&8&2&22& 26&26&20&0&20&
76&76&178\\ \hline $M_8(n)$&0&0&0&0& 2&6&0&0&0& 8&18&28&38&28&
0&0&0\\ \hline
\end{tabular}\] \pmn \begin{center}{Table 4}\end{center}
\subsection{$v=(\sqrt 2,\sqrt[3]{3})$}
\[\begin{tabular}{|c|c|c|c|c|c|c|c|c|c|c|c|c|c|c|c|c|c|}\hline
n&20&30&40&50&60&70&80&90&100&110&120&130&140&150&200&500&1000\\
\hline $S(n)$&4&6&6&7& 7&4&6&6&6& 7&7&7&5&7& 7&8&7\\ \hline
$M_4(n)$&0&0&0&0& 0&0&0&0&0& 0&0&0&0&4& 0&0&0\\ \hline
$M_5(n)$&0&4&12&14& 14&0&20&40&60& 60&40&20&0&6& 52&208&140\\
\hline $M_6(n)$&20&22&16&22& 38&70&48&28&8& 20&60&100&140&126&
96&84&720\\ \hline $M_7(n)$&0&4&12&14& 2&0&4&4&4& 0&0&0&0&14&
52&208&140\\ \hline $M_8(n)$&0&0&0&0& 6&0&8&18&28& 30&20&10&0&0&
0&0&0\\ \hline
\end{tabular}\] \pmn \begin{center}{Table 5}\end{center}
\subsection{$v=(\sqrt 2,e)$}
\[\begin{tabular}{|c|c|c|c|c|c|c|c|c|c|c|c|c|c|c|c|c|c|}\hline
n&20&30&40&50&60&70&80&90&100&110&120&130&140&150&200&500&1000\\
\hline $S(n)$&6&8&7&6& 6&8&7&6&6& 6&6&6&6&6& 6&8&7\\ \hline
$M_4(n)$&0&4&0&0& 0&0&0&0&0& 0&0&0&0&0& 0&0&0\\ \hline
$M_5(n)$&6&6&6&8& 14&14&4&14&14& 14&14&4&14&14& 38&92&52\\ \hline
$M_6(n)$&8&6&28&34& 32&42&72&62&72& 82&92&122&112&122&
124&316&896\\ \hline $M_7(n)$&6&14&6&8& 14&14&4&14&14&
14&14&4&14&14& 38&92&52\\ \hline
\end{tabular}\] \pmn \begin{center}{Table 6}\end{center}
\subsection{$v=(\sqrt[3]{2},e)$}
\[\begin{tabular}{|c|c|c|c|c|c|c|c|c|c|c|}\hline
n&50&100&110&120&130&140&150&200&500&1000\\ \hline $S(n)$&4&6&
7&7&7&7&7& 6&7&6\\ \hline $M_5(n)$&8&8& 8&8&8&8&8& 38&8&362\\
\hline $M_6(n)$&34&84& 94&104&114&124&134& 124&484&284\\ \hline
$M_7(n)$&8&8& 8&8&8&8&8& 38&8&346\\ \hline $M_8(n)$&0&0&
0&0&0&0&0& 0&0&8\\ \hline
\end{tabular}\] \pmn \begin{center}{Table 7}\end{center}
\subsection{$v=(e,\pi)$}
\[\begin{tabular}{|c|c|c|c|c|c|c|c|c|c|c|}\hline
n&50&100&110&120&130&140&150&200&500&1000\\ \hline
$S(n)$&4&12&13&11&10&12&13&10&11&6\\ \hline
$M_4(n)$&0&42&32&22&26&16&6&0&0&0\\ \hline
$M_5(n)$&6&22&22&22&4&40&80&56&58&142\\ \hline
$M_6(n)$&38&22&42&62&76&56&36&122&428&716\\ \hline
$M_7(n)$&6&0&0&0&10&14&14&8&0&142\\ \hline
$M_9(n)$&0&0&0&0&10&0&0&8&0&0\\ \hline
$M_{10}(n)$&0&0&0&4&4&12&0&6&12&0\\ \hline
$M_{11}(n)$&0&0&0&10&0&2&6&0&2&0\\ \hline
$M_{12}(n)$&0&0&12&0&0&0&8&0&0&0\\ \hline
$M_{13}(n)$&0&6&2&0&0&0&0&0&0&0\\ \hline
$M_{14}(n)$&0&8&0&0&0&0&0&0&0&0\\ \hline
\end{tabular}\] \pmn \begin{center}{Table 8}\end{center}
\section{Cautious conclusions}
\pn We see that numbers $S(n)$ demonstrate stabilization or a very
slow growth. It will be interesting to find out what really
happens. Also it will be interesting to perform analogous
computations in the three dimensional case. Authors intend to do
it. \pn \begin{rem} There exists another approach to a multidimensional
three gap theorem. See \cite{HM} and bibliography there. \end{rem}

\vspace{5mm}
\end{document}